\documentclass[english]{amsart}

\usepackage{babel}
\usepackage{amstext}
\usepackage{amsmath}
\usepackage{amsfonts}
\usepackage{latexsym}
\usepackage{ifthen}
\usepackage{xypic}
\xyoption{all}
\newcommand{\codim}{{\rm codim}}

\newcommand\sE{{\mathcal E}}

\newcommand\sF{{\mathcal F}}

\newcommand\sI{{\mathcal I}}

\newcommand\sL{{\mathcal L}}
\newcommand\sO{{\mathcal O}}

\newcommand\sC{{\mathcal C}}

\newcommand\bR{{\mathbb R}}
\newcommand\bZ{{\mathbb Z}}
\newcommand\bC{{\mathbb C}}
\newcommand\bQ{{\mathbb Q}}

\newcommand\sB{{\mathcal B}}

\newcommand\bP{{\mathbb P}}

\newcounter{lemma}

\newtheorem{lemma1}[lemma]{\setcounter{equation}{0}}

\newenvironment{lemma}{\begin{lemma1}{\bf Lemma.}}{\end{lemma1}}

\newenvironment{example}{\begin{lemma1}{\bf Example.}\rm}{\end{lemma1}}

\newenvironment{theorem}{\begin{lemma1}{\bf Theorem.}}{\end{lemma1}}
\newenvironment{question}{\begin{lemma1}{\bf Question.}}{\end{lemma1}}

\newenvironment{proposition}{\begin{lemma1}{\bf Proposition.}}{\end{lemma1}}

\newenvironment{corollary}{\begin{lemma1}{\bf Corollary.}}{\end{lemma1}}

\newenvironment{remark}{\begin{lemma1}{\bf Remark.}\rm}{\end{lemma1}}
\newenvironment{definition}{\begin{lemma1}{\bf Definition.}}{\end{lemma1}}

\newenvironment{notation}{\begin{lemma1}{\bf Notation.}}{\end{lemma1}}

\newenvironment{conjecture}{\begin{lemma1}{\bf Conjecture.}}{\end{lemma1}}

\newenvironment{Induction Step}{\begin{lemma1}{\bf Induction Step.}}{\end{lemma1}}
\newenvironment{Proof of Theorem 1.2}{\begin{lemma1}{\bf Proof of Theorem 1.2.}}{\end{lemma1}}

\title {Submanifolds with ample normal bundles and a conjecture of Hartshorne} \author{Thomas Peternell}
\date{\today}

\begin{document}

\maketitle

\tableofcontents

\section{Introduction} We consider a projective manifold $Z$ and submanifolds $X$ and $Y$ with {\it ample normal bundles. }
In [Ha70,chap III,4.5] R.Hartshorne stated the following conjecture:

\begin{conjecture} If $\dim X + \dim Y \geq \dim Z$, then $X \cap Y \ne \emptyset.$ 
\end{conjecture} 

Although there are some remarkable positive results, the conjecture is in principle wide open; see section 2 for a short
description of what is known so far.\\ 
We observe  first in this paper that the conjecture holds generically. To be more precise we introduce 
$$ VG_k(Z) \subset Z $$
to be the set of points $x$ such given an irreducible $k-$cycle through $x,$ then a multiple of the cycle moves in a family covering $Z.$ 
From general properties of the Chow scheme, it is clear that $Z \setminus VG_k(Z)$ is a countable union of proper subvarieties of $Z.$ 
Using criteria of Barlet resp. Fulton-Lazarsfeld to decide when $X$ and $Y$ meet, one deduces the following

\begin{theorem} Suppose that under the above conditions that $X \cap VG_{m-1}(Z) \ne \emptyset, $ where $ m = \dim X$ and that $N_Y$ is even positive in the sense
of Griffiths. Then $X \cap Y \ne \emptyset.$ If we make the stronger assumption that  $X \cap VG_{m}(Z) \ne \emptyset, $ then the ampleness assumption 
on $N_X$ can be dropped. 
\end{theorem}

In particular, if $X$ contains a sufficiently general point of $Z$, then the Hartshorne conjecture holds for $X$ and any $Y$. 
The a priori stronger condition that $N_Y$ is
positive in the sense of Griffiths is needed to ensure that $Z \setminus Y$ is $k-$convex (in the sense of Andreotti-Grauert), where $k = \codim Y.$ 
Vector bundles which are positive in the sense of Griffiths are necessarily ample, but it is still unknown whether the converse is 
also true. \\
Instead of $VG(Z)$ we can also use either $VG^a_k(Z)$ or $VG^{sm,a}_k(Z)$, using only $k-$cycles with ample normal bundles resp. smooth 
$k-$cycles with ample normal bundles and everything said so far remains true. 

\vskip .2cm \noindent We next observe that the Hartshorne conjecture holds if the class $[X]$ is contained in the interior of the cone generated by the
cohomology classes of irreducible $m-$dimensional subvarieties, where $m = \dim X.$ We discuss this property in detail in section 4. 
In particular we show that the Hartshorne conjecture holds once the following question has a positive answer: 
\vskip .2cm \noindent
{\it given a submanifold $X_m \subset Z$ with ample normal bundle, is the class $[X]$ an interior
point of the cone generated by the cohomology classes of $m-$dimensional subvarieties of $Z?$}
\vskip .2cm \noindent 
As to positive results, we verify the Hartshorne conjecture in the following cases - we always assume $N_Y$ to be positive in the sense of Griffiths. 
\vskip .2cm \noindent
\begin{itemize} 
\item $Z$ is a degree 2 cover over a projective homogeneous manifold;
\item $Z$ is a smooth hypersurface in a weighted projective space;
\item $Z$ admits a fibration over a curve whose general fiber is homogeneous;
\item $Z$ is a $\bP_1-$bundle over a threefold;
\item $Z$ is a $4-$fold and $X$ is a surface with $\kappa (X) = - \infty$, resp. 
\item $Z$ is a $4-$fold and $X$ is a non-minimal surface; moreover $\kappa(\sO_Z(D)) \geq 1$ for every effective divisor $D$;
\item $Z$ is a Fano manifold of index $n-1,$ i.e., $Z$ is a del Pezzo manifold;
\item $Z$ is a Fano manifold of index $n-2$ with a few possible exception (only one exception in all dimensions $\geq 5$
, namely the intersection of three quadrics
in $\bP_{n+3}).$
\end{itemize}

\section{Update on the Hartshorne Conjecture}
\setcounter{lemma}{0}

In this section we collect the known results on the Hartshorne conjecture and fix the following situation: 
\vskip .2cm \noindent 
{\it $Z$ is a projective manifold of dimension $n$ with compact submanifolds $X,Y$ of dimensions
$m,k$ such that $m+k \geq n.$ We assume that the normal bundle $N_X$ is ample and $N_Y$ is Griffiths-positive.}
\vskip .2cm \noindent
Recall that a rank $r-$vector bundle $E$ is said to be positive in the sense of Griffiths, {\it G-positive} for short,
if there is a hermitian metric  on $E$ such that the
curvature $\Theta$ of the canonical connection fulfilles the following positivity condition
$$ \sum_{i,j\alpha,\beta} \Theta_{i,j,\alpha,\beta}(z) \zeta^{\alpha}\overline {\zeta^{\beta}} \eta^{i} \overline {\eta^{j}} > 0$$
for all $z \in Z, (\zeta^{\alpha}) \in \bC^r \setminus \{0\}, (\eta^{i}) \in \bC^n \setminus \{0\}. $ 
\vskip .2cm \noindent
Notice first that if $m = n-1,$ then the conjecture is obviously true so that usually we shall assume $m \leq n-2.$ 
In [Ba87] and [BDM94] the most general result is proved - generalizing [Lu80] and [FL82].

\begin{theorem} If $Z$ is a hypersurface in a homogeneous manifold, then $X \cap Y \ne \emptyset. $ 
(It suffices that there is an open neighborhood of $X$ in $Z$ which is biholomorphic to a locally closed hypersurface 
of a homogeneous manifold). 
If $Z$ is a hypersurface in $\bP_{n+1},$ then it suffices both normal bundles to be ample. 
\end{theorem} 

The paper [DPS90] deals with special 4-folds:
\begin{theorem} If $Z$ is a $\bP_2-$bundle over a smooth projective surface, then $X \cap Y \ne \emptyset.$ 
\end{theorem} 

In the thesis [Poe92] the Hartshorne conjecture (for G-positive normal bundles) is settled for many $\bP_1-$bundles over
threefolds - the general case will be done in section 5.
The last result is due to Migliorini [Mi92].

\begin{theorem} Assume that $\dim Z = 4$ and that $b_2(Z) = 1$ or that $Z$ is a complete intersection in some projective space. 
The normal bundles
$N_X$ and $N_Y$ are supposed only to be ample. Suppose furthermore that the surface $X$ is minimal of non-negative Kodaira dimension and
that $c_1(N_X)^2 > 2 c_2(X).$ Then $X \cap Y \ne \emptyset.$
\end{theorem} 

A  general notice: by taking hyperplanes section or submanifolds in $X$ or $Y$ with ample normal bundles, we can always reduce - and do - to the case
$$ \dim Z = \dim X + \dim Y. $$

\section{The generic Hartshorne Conjecture}
\setcounter{lemma}{0}

We fix again a projective manifold $Z$ of dimension $n$ and submanifolds $X$ and $Y$ with $\dim X = m$ and $\dim Y = k$ subject to the 
condition $m + k = n$.

Here is a criterion due to Barlet to verify $X \cap Y \ne \emptyset.$

\begin{proposition} Assume that there is an effective divisor ($ = (m-1)-$cycle) $D \subset X$ moving in an irreducible family $(D_t)$ of $(m-1)-$cycles in $Z$ such that  
$D_{t_0} \cap Y \ne 0$ for some $t_0.$ If $N_X$ is ample and $N_Y$ is G-positive,  then $X \cap Y \ne \emptyset.$
\end{proposition} 

For the proof we refer to [Ba87], [BPS90,1.4,1.5] for the case that $N_X$ and $N_Y$ both G-positive,
and to [BDM94], [Ba99] in the case that $N_X$ is merely ample (and $N_Y$ G-positive).  

The other - related - criterion is due to Fulton-Lazarsfeld [FL82], [Fu84,12.2.4]:  

\begin{proposition} Suppose $N_Y$ ample and $m + k = n.$ 
Suppose furthermore that  $X$ is numerically equivalent to an effective cycle meeting $Y$, then 
$X \cap Y \ne \emptyset$ (here we do not assume $N_X$ to be ample). 
\end{proposition} 

Barlet's criterion has the advantage that one needs ``only'' to move divisors in $X$, on the other hand the assumptions are stronger.     

\begin{notation} {\rm We denote the cycle space of $Z$ by $\sB(Z)$ and by $\sB_k(Z)$ the subspace of $k-$cycles. If $S \subset \sB_k(Z)$ is an 
irreducible subvariety, we consider the
associated family $q: \sC_S \to S$ with projection $p: \sC_S \to Z.$ If $D \subset Z$ is a $k-$cycle, we consider the associated 
point $[D] \in \sB(Z)$ and a positive-dimensional irreducible subvariety $S \subset \sB(Z)$ (usually an irreducible component) containing $[D]$. 
We say that $D$ deforms in the family
$\sC_S$, or, introducing $D_s = p(q^{-1}(s))$ (as cycle), that $D$ deforms in the family $(D_s).$ This family is {\it covering} if $p$ is
surjective. }
\end{notation} 

Following Koll\'ar [Ko95] we define - however in a somehow different setting - very general points. 

\begin{definition} $VG_k(Z)$ is the set of points $z \in Z$ subject to the following condition.
If $D$ is any irreducible $k-$dimensional subvariety passing through $z$, then some multiple $mD$ moves in a family covering $Z.$ 
\end{definition} 

Similarly as in [Ko95] we have

\begin{proposition} For all $k > 0$ there are at most countably many irreducible subvarieties $W_r \subset Z$ such that $Z \setminus VG_k(Z) \subset
\bigcup_j W_j. $
\end{proposition} 

\begin{proof} Let $S_j \subset \sB_k(Z)$ denote those irreducible components for which the projection $p_j: \sC_{S_j} \to Z$ is
not surjective. These are at most countable many, simply because $\sB_k(Z)$ has only countably many components. Now set
$$ W_j = p_j(\sC_{S_j}). $$ 
So if $x \in Z \setminus \bigcup_j W_j$ and if $D$ is an irreducible subvariety containing $x$, then for any component $S$ of $\sB_k(Z)$
containing $[D]$, then $S \ne S_j$ for all $j$ and therefore the associated family covers $Z$ so that even $D$ moves in a covering
family.
\end{proof}

Putting things together we obtain

\begin{theorem} 
\begin{enumerate} 
\item Suppose in our setting that $N_X$ is ample and that $N_Y$ is G-positive. If $X \cap VG_{m-1}(Z) \ne \emptyset,$ then $X \cap Y \ne \emptyset.$ 
\item If $N_Y$ is merely ample (without any assumption on $N_X$), and if $X \cap VG_m(Z) \ne \emptyset,$ then $X \cap Y \ne \emptyset. $ 
\end{enumerate}
In particular there is a countable union $T$ of subvarieties of $Z$ having the following property. If $X$ and $Y$ are
submanifolds of $Z$ with ample normal bundles and $\dim X + \dim Y \geq \dim Z$ such that $X \not \subset T$, then 
$X \cap Y \ne \emptyset.$ 
\end{theorem} 

\begin{proof} (1) By (3.1) we need to move some irreducible divisor $D \subset X$ to meet $Y.$ Choose $x \in X \cap VG_{m-1}(Z)$ and
take any irreducible divisor $D \subset X$ passing through $x.$ Then $D$ moves in a family covering $x,$ hence some deformation
of $D$ meets $Y$ and we conclude.  \\
(2) Choose $x \in X \cap VG_m(Z).$ Then $X$ moves in a family covering $Z.$ Now apply (3.2) to conclude. 
\end{proof} 

It is actually not necessary to work with singular cycles; we can define $VG_k^{sm}(Z)$ as the set of points $z \in Z$ with the property that if $D$ is
a $k-$dimensional smooth subvariety passing through $z,$ then some multiple of $D$ moves in a family covering $Z$. 
Then all what we said for $VG_k(Z)$ remains true for $VG_k^{sm}(Z).$  
We can even put more conditions on the cycles, namely we can ask $D$ to have ample normal bundle (or rather ample normal sheaf) in $Z$. The resulting sets are denoted
$VG_k^a(Z)$ resp. $VG_k^{sm,a}(Z).$ 
 
\vskip .2cm \noindent In general it is difficult to compute $VG_k(Z),$ even in the simplest case $\dim Z = 2$ and $k = 1.$ So suppose $Z$ a
projective surface and suppose $VG_1(Z) = Z.$ Then $Z$ does not contain any irreducible curve $C$ with $C^2 < 0,$ in particular $Z$ is minimal.
Moreover: 
\begin{itemize}
\item $\kappa(Z) = - \infty$ iff $Z = \bP_2, \bP_1 \times \bP_1$ or $Z = \bP(E)$ with $E$ a semi-stable rank $2-$bundle over a curve $B$ of genus $\geq 2$ or 
of the form $\sO_B \oplus L$ with $L$ torsion; 
\item $\kappa  (Z) = 0$ iff $X$ is torus, hyperelliptic or K3/Enriques without $(-2)-$curves. 
\end{itemize}
\vskip .2cm \noindent 
If however we consider $VG_1^a(Z)$, things gets much easier: obviously
$$ VG^a_1(Z) = Z $$
for all surfaces $Z.$ 
At the moment I do not have any example of a threefold or a fourfold $Z$ such that $VG^a_1(Z) \ne Z.$

\section{Some general observations}
\setcounter{lemma}{0}

\begin{notation} {\rm Let $Z_n$ be a projective manifold. Then $K^a_r(Z)$ denotes the closed cone of classes of effective $r-$cycles $\sum a_i W_i$ 
(with $W_i$ irreducible of dimension $r$) in $A_r(Z),$ in the Chow ring of $Z$.\\ 
If we consider numerical instead of rational equivalence, we obtain the cone
$K_r(Z) \subset H^{n-r,n-r}_{\bf R}(Z).$ The class numerical of $W_i$ will be denoted by
$$ [W_i] \in  H^{n-r,n-r}_{\bf R}(Z). $$
Given subvarieties $X$ and $Y$ such that $\dim X + \dim Y = n,$ we can form the intersection product 
$$ X \cdot Y \in A_0(X) \simeq \bZ $$
which will always be considered as a number. 
}
\end{notation} 

\begin{theorem} Let $X$ and $Y $be submanifolds of $Z$ of dimensions $m$ and $k$ with $m + k = \dim Z = n$.
\begin{enumerate} 
\item If $N_X$ or $N_Y$ is ample, then $X \cdot Y = 0$ if and only if $X \cap Y = \emptyset.$ 
\item If $N_Y$ is ample and
if $X \cap Y = \emptyset,$ then $[X] \in \partial K^a_n(Z)$. 
\end{enumerate} 
\end{theorem}

\begin{proof} 
(1) If $X \cap Y = \emptyset,$ then of course $X \cdot Y = 0.$ The other direction is [FL82, Theorem 1]. \\
(2) We consider the linear form 
$$ \Phi_Y: A_m(Z) \to \bZ,$$
$$ \sum a_i [W_i] \mapsto \deg (\sum a_i (Y \cdot W_i))$$
(where $Y \cdot W_i \in A_0(W_i) \simeq \bZ$).
By [FL82], the ampleness of $N_Y$ implies that $\Phi_Y(W) \geq 0 $ for $W$  irreducible of dimension $m.$ Thus 
$\Phi_Y \vert K^a_m(Z) \geq 0.$ Now
$$ \Phi_Y(X) = X \cdot Y = 0.$$
Thus $[X]$ cannot be in the interior of $K^a_n(Z),$ since $\Phi \ne 0.$ 
\end{proof} 

\begin{corollary}  Let $X$ and $Y $be submanifolds of $Z$ of dimensions $m$ and $k$ with $m+k = \dim Z.$ Suppose $N_X$ and $N_Y$ ample. 
If $X \cap Y = \emptyset,$ then $[X] \in \partial K_m(Z)$ and $[Y] \in \partial K_k(Z).$ 
\end{corollary} 

\begin{proof} We just have to notice that for $W_1, W_2$ numerically equivalent, we have $\deg (Y \cdot W_1) = \deg (Y \cdot W_2) $.
Indeed, in $H^*(Z,\bR) $ one has $\deg (Y \cdot W_j) = [Y] \cdot [W_j],$ see [Fu84,chap.19]. 
\end{proof} 

The Hartshorne conjecture would therefore be a consequence of a positive answer to the following question. 

\begin{question} Let $X \subset Z$ be a submanifold of dimension $m$ in the projective manifold $Z.$ If $N_X$ is ample, must $[X]$ be in the interior of $K_m(Z)?$
\end{question} 

In codimension 1, the answer is easy, since a big divisor is the sum of an ample and an effective $\bQ-$divisor: 

\begin{proposition} Assume $X \subset Z_n$ is a smooth divisor with ample normal bundle. Then $[X] \in K_{n-1}(Z)^{o},$ the interior of the
pseudo-effective cone of $Z.$
\end{proposition} 

Even in dimension 1, the analogous statement is open:
let $X \subset Z_n$ be a smooth compact curve with ample normal bundle. Is  
$$[X] \in K_1(Z)^{0} = NE(Z)^{0}?$$
This comes down to solve positively the following problem:
\vskip .2cm \noindent
{\it Let $L$ be a nef line bundle and a smooth curve $C \subset Z$ with ample normal bundle. If $L \cdot C = 0,$ is $L \equiv 0?$ }
\vskip .2cm \noindent 
If $\dim Z = 2,$ this follows from Hodge Index Theorem. Here are some partial results in higher dimensions. 

\begin{proposition} Let $Z$ be a projective manifold, $L $ a nef line bundle on $Z$ and $C \subset Z$ a smooth curve with ample normal bundle.
If $L \cdot C = 0,$ then $\kappa (L) \leq 0.$
\end{proposition} 

\begin{proof} By  [PSS99], there is a positive number $c$ such that for all $t:$
$$ h^0(tL) \leq \sum_{k=0}^cth^0(S^kN^*_C \otimes tL_C). $$ 
Since $L \vert C \equiv 0,$ we obtain 
$$ h^0(tL) \leq h^0(tL_C) \leq 1. $$
Hence $\kappa (L) \leq 0.$
\end{proof} 

If $Z$ is a $\bP_{n-1}-$bundle over a curve, things are easy:

\begin{proposition} Let $p: Z = \bP(\sE) \to B$ be a $\bP_{n-1}-$bundle over the smooth compact curve $B.$ Let $C \subset Z$
be a smooth curve with ample normal bundle. Then $[C] \in K_1(Z)^{0}.$
\end{proposition} 

\begin{proof} By possibly passing to a covering of $B,$ we may assume that $C$ is a section of $p.$ Then $C$ corresponds
to an epimorphism
$$ \sE \to \sL \to 0 $$
(such that $C = \bP(\sL)$). Let $\sF$ denote its kernel. Then 
$$ N_{C/Z} \simeq \sL \otimes \sF^*,$$
hence $\sL \otimes \sF^*$ is ample. In order to prove our assertion, we pick a nef $\bQ-$divisor $D$ and need to show that
$D \cdot C > 0.$ 
We may write $$ D \equiv \zeta + p^*(A), $$
where $\zeta = \sO_{\bP}(\sE)$ and $A$ is a $\bQ-$divisor on $B.$ 
From the exact sequence
$$ 0 \to \sF \otimes \sL^* \to \sE \otimes \sL^* \to \sO_Z \to 0$$
and the ampleness of $\sL \otimes \sF^*$, we deduce $\deg A > \deg \sL^*$ (notice that $\sE \otimes \sL^*$ cannot be nef). 
Thus
$$ D \cdot C = \zeta \cdot C + \deg A = \deg \sL + \deg A > 0.$$
\end{proof} 

\begin{theorem} Let $C \subset Z$ be a smooth curve with ample normal bundle. Suppose $C$ moves in a family $(C_t)$ which covers $Z.$ 
Then $[C]$ is in the interior of $\overline{NE}(Z).$ 
\end{theorem} 

\begin{proof} We must show that, given a nef line bundle on $Z$ with $L \cdot C = 0,$ then $L \equiv 0.$ 
Consider the nef reduction $f: Z \dasharrow S$ of $L$, see [workshop]. Thus $f$ has the following properties.
\begin{itemize}
\item $f$ is almost holomorphic, i.e. the general fiber $F$ of $f$ is compact;
\item $L \vert F \equiv 0$;
\item If $B$ is any curve through a general point of $Z$, then $\dim f(B) = 0 $ iff $L \cdot B = 0.$ 
\end{itemize} 
Since $N_C$ is ample, so does $N_{C_t}$ for general $t.$ On the other hand $\dim f(C_t) = 0,$ since $L \cdot C_t = 0.$ Both facts 
together are in
contradiction unless $S$ is a point. But then $L \equiv 0$. 
\end{proof} 

This result remains true for singular curves assuming that the normal sheaf $(\sI_C/\sI_C^2)^*)$ is ample. Even if $\kappa (L) = 0,$ the general problem however is open; 
specifically we ask

\begin{question} Let $Z$ be a smooth projective threefold, $Y \subset Z$ a smooth hypersurface with nef normal bundle and $C \subset Z$ a smooth curve
with ample normal bundle. Is $Y \cap C \ne \emptyset ?$
\end{question}

\begin{example} {\rm In [FL82] Fulton and Lazarsfeld gave an example of a surface $X$ in a $4-$fold $Z$ with ample normal bundle such that no multiple of 
$X$ moves inside $Z.$ Here we show that nevertheless $[X]$ is in the interior of $K_2(Z).$ \\ 
Let $\sF$ be an ample rank 2-vector bundle on $\bP_2$ given by an exact sequence 
$$ 0 \to \sO(-n)^2 \to \sO(-1)^4 \to \sF \to 0 $$
for a suitable large $n.$ These bundles were constructed by Gieseker [Gi71]. We consider
$$ Z = \bP(\sO \oplus \sF^*) $$
with projection $\pi: Z \to \bP_2$
and the submanifold
$$ X = \bP(\sO) \simeq \bP_2.$$
Note that the normal bundle $N_{X/Z} \simeq \sF$ is ample. In [FL82] it is shown that no multiple of $X$ moves in $Z.$ 
Consider a line $l \subset X.$ Then the normal bundle $N_{l/Z}$ is ample and the deformations of $l$ cover $Z.$ Hence
by (3.1) $X$ meets every surface $Y \subset Z$ with G-positive normal bundle. \\
We prove that Question 4.4 has a positive answer for $X$:
$$ [X] \in K_2(Z)^{o}.$$
Consider now a {\it general} line $l \subset \bP_2.$ Since $\sF$ is stable (this is obvious from $H^0(\sF) = 0$), the Grauert-M\"ulich
theorem determines the splitting behaviour:
$$ \sF \vert l = \sO(n-2) \oplus \sO(n-2).$$ Therefore
$$ Z_l := \pi^{-1}(l) \simeq \bP(\sO \oplus \sO(2-n)^2). $$
Consider the map
$$ \phi: H^4(Z,\bR) \to H^4(Z_l,\bR) $$
given by $S \mapsto S \cap Z_l.$ Then $\phi(K_2(Z)) \subset  \overline{NE}(Z_l). $
Let $K' = \phi(K_2(Z))$, a closed subcone of the 2-dimensional cone $\overline {NE}(Z_l).$ 
It is immediately seen that one of the two boundary rays of $K'$ is occupied by a line $l'$ in a fiber of $Z_l \to l.$ 
Set
$$ X_l = X \cap Z_l.$$ 
This is the section $\bP(\sO_l) \subset Z_l$ and it has normal bundle 
$$N_{X_l/Z_l} = \sO(n-2) \oplus \sO(n-2).$$ 
Let $\zeta = \bP(\sF). $ Since $\zeta$ is ample, we find $m > 0$ and an element
$$ S \in \vert m\zeta \vert, $$
where $p: \bP(\sF) \to \bP_2$ is the projection. 
We have an embedding $$S \subset \bP(\sF) \simeq \bP(\sF^*) \subset Z. $$
Let $S_l = S \cap Z_l,$ a multisection of $Z_l$ which is disjoint from $X_l.$ 
Since $l'$ and $S_l$ are independent in $H^4(Z_l,\bR) $, we can write 
$$ X_l = \mu S_l + \nu l' \eqno(*) $$
in $H^4(Z_l)$ with real (actually rational) coefficients $\mu, \nu.$ 
We claim that $\mu,\nu > 0,$ so that $[X_l ]$ is not extremal in $K'.$ Hence $[X]$ is not extremal in
$K_2(Z),$ i.e., contained in the interior of $K_2(Z). $ 
To verify the positivity of $\mu$ and $\nu$ we first dot (*) with a $\pi-$ fiber $F$ to obtain $\mu = {{1} \over {m}}. $
Then we dot with $\bP(\sF^*_l) $ and use 
$$ S_l \cdot \bP(\sF^*_l) = - \zeta_{\sF_l}^2 < 0 $$
to deduce
$$ \nu = - {{1} \over {m}}  S_l \cdot \bP(\sF^*_l) > 0.$$
Thus we conclude that $[X] \in K_2(Z)^{o}.$
} 
\qed
\end{example}

We next prove a statement which would be an immediate consequence of a positive answer to the Hartshorne conjecture. 
\begin{theorem} Let $X,Y \subset Z$ be  compact submanifolds of dimensions $m$ and $n.$ Assume $\dim Z = m+n.$ Suppose that $X$ and $Y$ meet 
transversally in $d$ points $x_1, \ldots, x_d.$ Let $\pi: \hat X \to X$ be the blow-up of $x_1, \ldots, x_d$ with exceptional divisors $E_j.$ Let $\hat X$ and $\hat Y$ be
the strict transform of $X$ and $Y$. Then at least one of the normal bundles $N_{\hat X}$,$N_{\hat Y}$ is not G-positive. 
Hence $\pi^*(N_X) \otimes \sO_{\hat X}(- \sum E_j \vert \hat X)$ or  $\pi^*(N_X) \otimes \sO_{\hat Y}(- \sum E_j \vert \hat Y)$ is not
G-positive. 
\end{theorem} 

\begin{proof} We argue by contradiction and need to construct a divisor $D \subset \hat X$ which moves in a family $(D_t)$ such 
that $D_{t_0} \cap \hat Y $ for some $t_0.$ 
We consider the exceptional divisor $E_1$ lying over $x_1$ and put $D = E_1 \cap \hat X.$ Then $D \simeq \bP_{m-1}$ is a linear subspace, 
and since $E_1 \cap \hat Y \ne \emptyset,$ some deformation  of $D$ in $E$ meets $\hat Y.$ Hence not both $N_{\hat X}$ and $N_{\hat Y}$ can be
Griffiths-positive by (3.1). 

\end{proof} 

For later use we establish the Hartshorne conjecture for degree 2 covers of homogeneous manifolds.

\begin{theorem} Let $Z$ be a projective manifold with a degree 2 cover $f: Z \to W$ over a projective homogeneous manifold $W.$ 
Let $X,Y \subset Z$ submanifolds with $N_X $ ample, $N_Y$ G-positive and $\dim X + \dim Y \geq \dim Z.$ 
Then $X \cap Y \ne \emptyset.$ 
\end{theorem} 

\begin{proof} Let $Y' = f(Y).$ Since $W$ is homogeneous, $Y'$ moves in a family covering $W.$ Hence $f^*(Y') = f^{-1}(Y')$, the 
scheme-theoretic preimage,
moves in a family $(f^*(Y'_t))$ covering $Z.$ 
Thus for some $t,$ we have $X \cap f^{-1}(Y'_t) \ne \emptyset.$ From (3.2) it follows
$$ X \cdot f^{*}(Y') = X \cdot f^{*}(Y'_t)) > 0.$$
If $\deg f \vert Y = 2,$ or if $Y$ lies in the branch locus $B$ of $f,$ 
then $f^{-1}(Y') = Y$ set-theoretically, hence $X \cdot Y > 0,$ so that $X \cap Y \ne \emptyset.$ If $\deg f \vert Y = 1$ and 
if $Y \not \subset B,$ 
then $f^{*}(Y')$ has
a second component $\tilde Y.$ Assume $X \cap Y = \emptyset.$ Then $X \cdot \tilde Y > 0,$ so that $X \cap \tilde Y \ne \emptyset.$ 
We now show that there is a divisor $D \subset Y$ which is also contained in $\tilde Y$ deforming in a covering family of $\tilde Y$. 
Some deformation will therefore meet $X$,
so that by (3.1) we arrive at a contradiction. In order to produce $D$, we consider the ramification divisor $R \subset W.$ 
Since $W$ is homogeneous, $R$
moves in a covering family. Hence $R \cap Y'$ moves in a family $(D'_t)$ covering $Y'$ with $D_0 =  R \cap Y'.$ 
Now consider the family $f^{*}(D_t)$ in 
$Y \cup \tilde Y$; notice $f^*(D_0) \subset Y \cap \tilde Y.$ Furthermore for general $t$ we can write 
$$ f^*(D_t) = S_t \cup \tilde S_t$$
with $S_t \subset Y_t$ and $\tilde S_t \subset \tilde Y.$ 
The family $(f^{*}(D_t) \cap \tilde Y)$ thus deforms a divisor 
contained in $\tilde Y \cap Y$, namely ${{1} \over {2}}f^*(D_0)$ to a divisor in $\tilde Y$ meeting $X$, and we are done.  

\end{proof} 

\begin{theorem} Let $Z$ be a projective manifold of dimension $n$ and $f: Z \to B$ be a surjective map with connected fibers to a 
smooth curve $B.$ 
Assume that the general fiber $F$ of $f$ is homogeneous. 
Let $X,Y \subset Z$ be submanifolds with $N_X  $ ample, $N_Y$ G-positive and $\dim X + \dim Y = n.$ Then $X \cap Y \ne \emptyset.$
\end{theorem} 

\begin{proof} By the ampleness of $N_X$ amd $N_Y$ the maps $f \vert X$ and $f \vert Y$ are onto $B.$ Thus $F \cap X$ and $F \cap Y$ 
are divisors in 
$X$ resp. $Y.$ We want to move $F \cap X$ inside $F$ to meet $F \cap Y$. But this is obvious by homogeneity. Now we conclude by (3.1). 
\end{proof} 

This theorem applies e.g. to manifolds $Z$ with $\kappa (Z) = 1 $ such that its Iitaka fibration is holomorphic with general fiber
a torus. 
\section{Fourfolds and Fano manifolds} 
\setcounter{lemma}{0}

We first show that the Hartshorne conjecture holds for $\bP_1-$bundles over threefolds.

\begin{theorem} Let $Z$ be a smooth projective $4-$fold, $\pi: Z \to W$ a $\bP_1-$bundle. Let $X \subset Z$ and $Y \subset Z$ be
surfaces with G-positive normal bundles. Then $X \cap Y \ne \emptyset.$ 
\end{theorem}

\begin{proof} After a finite \'etale cover of $W$ we may write
$$ Z = \bP(E) $$
with a rank $2-$bundle $E$ on $W.$  Passing to $\bQ-$bundles $E$, we may also assume 
$$ c_1(E) = 0.$$
So from now on, all bundles are $\bQ-$bundles. 
It is easy to see ([Poe92]), that $\pi \vert X$ and $\pi \vert Y$ are finite and that $X' = \pi(X)$ and $Y' = \pi(Y)$ are surfaces
with ample normal bundles in $W$. Thus $X'$ and $Y'$ meet in finitely many curves $C_j.$ \\
Let $$ \zeta = \sO_{\bP(\sE)}(1). $$
The equation $c_1(E) = 0$ implies via the Hirsch-Leray relation $\zeta^2 = -\pi^*(c_2(E)). $ Therefore we may write in $N^*(Z): $
$$ X \equiv \zeta \cdot \pi^*(D) + \pi^*(C) \eqno (1) $$
with $D \in N^1_{\bQ}(W)$ and $C \in N^2_{\bQ}(W) \simeq N_1^{\bQ}(W). $ In other words $D \equiv \sum a_i D_i$ with irreducible 
hypersurfaces $D_i \subset W$; $a_i 
\in \bQ$
and $C \equiv \sum b_j C_j$ with irreducible curves $C_j \subset W$ and $b_j \in \bQ.$ 
\vskip .2cm \noindent We are going to fix some notation. We consider an irreducible, possibly singular, curve $C \subset W$ and the 
ruled surface
$Z'_C = \pi^{-1}(C)$ whose normalization is denoted by $\nu: Z_C \to Z'_C.$ Using the notations of [Ha77,V.2], the surface $Z_C$ has 
an invariant $e$ and a section $C_0$
of minimal self-intersection $C_0^2 = -e. $ We also have 
$$ \zeta = C_0 + {{e} \over {2}} F, $$
where $F$ is a ruling line.  
\vskip .2cm \noindent 
(A) Suppose that there is an ample line bundle $L$ on $W$ such that
$$ X \cdot \pi^*(L) \cdot \zeta > 0. \eqno (2) $$
We may assume $L$ very ample, take a general element $S \in \vert L \vert $ and put $C = S \cap X'.$ 
Let $X_C = \nu^{-1}(X \cap Z'_C)$ and $\zeta_C = \nu^*(\zeta \vert Z'_C).$ 
Writing $$ X_C = C_0 + \mu F,$$ 
equation (2) reads 
$$ X_C \cdot \zeta_C = \mu - {{e} \over {2}} > 0.$$
Using the description of the pseudo-effective and the nef cone of a ruled surface as give in [Ha77,V.2], we conclude that $X_C$ is a 
big divisor in
$Z_C.$ Therefore a multiple of $X_C$ moves to fill up $Z_C.$ Hence a multiple of $X \cap Z'_C$ moves and fills up $\pi^{-1}(X'),$ since 
we may 
also vary $C.$ Since $\pi^{-1}(X') \cap Y \ne \emptyset,$ we may apply Theorem 3.1 and conclude $X \cap Y  \ne \emptyset.$ 
\vskip .2cm \noindent
(B) So we may assume that
$$ X \cdot \pi^*(L) \cdot \zeta \leq 0 \eqno (3) $$
for all ample $L$ on $W.$ 
Putting (1) into (3) gives 
$$ L \cdot C \leq 0$$
for all ample $L$ on $W$. Thus $-C \in \overline{NE}(W).$
Using again (1),
$$ X^2 = (\zeta \cdot \pi^*(D) + \pi^*(C))^2 = 2 \zeta \cdot \pi^*(D) \cdot \pi^*(C) = 2 X \cdot \pi^*(C).$$
The ampleness of $N_X$ implies $X^2 > 0,$ hence $X \cdot C > 0.$ 
By the projection formula
$$ X \cdot \pi^*(C) = d X' \cdot C,$$
where $d$ is the degree of $X$ over $X'.$ Hence
$$ X' \cdot C > 0. $$
On the other hand, $-C \in \overline{NE}(W),$ which leads to a contradiction, the divisor $X'$ being nef in $W.$

\end{proof}

\begin{remark} Theorem 5.1 should of course also be true if the normal bundles are just ample. If $\deg \pi \vert X \geq 2$ and 
$\deg \pi \vert Y \geq 2$ and if every big and semi-ample divisor on $W$ is actually ample, 
this is seen as follows. We shall use the notations of the proof of (5.1) and
argue that if $\pi \vert X$ has degree at least $2$, then we have
$$ X \cdot \pi^*(L) \cdot \zeta \geq 0 \eqno (1)$$
for all ample line bundles $L$ on $W$. 
This is done using the computations in (5.1) by choosing a curve $C$ as intersection $S \cap X'$ with $S$ a general element in 
$\vert mL \vert.$ Then we use the theory of
ruled surfaces, applied to $Z_C$, to compute. \\
Next we claim that - assuming $X \cap Y = \emptyset$ - 
$$ X \cdot \pi^*(Y') \cdot \zeta = 0.  \eqno (2)$$
This is seen as follows. We take one of the irreducible curves $C_j \subset X' \cap Y'$ and form the ruled surface $Z_j = Z_{C_j}. $ 
Then $X_j $ and $Y_j$ are disjoint multi-sections - if we assume $X \cap Y = \emptyset$ - possibly reducible. By (3.1) no deformation 
of a multiple of any
component of $X_j$ 
meets $Y_j$ and vice versa. Using again [Ha77,V.2], this is only possible when $e = 0$ and $X_j,Y_j$ are sections with 
$X_j^2 = Y_j^2 = 0.$
This implies (2). \\
Now by our assumption the a priori only big and semi-ample divisor $Y'$ is ample. Therefore equation (1) and (2) together
yield
$$ X \cdot \pi^*(L) \cdot \zeta = 0 $$ 
for all ample line bundles $L$. Hence $L \cdot C = 0$ for all $L$ and therefore $C \equiv 0.$ 
Consequently $X^2 = \zeta^2 \cdot \pi^*(D^2) = 0,$ contradicting the ampleness of $N_X.$ 
\end{remark}

In the next theorem we put some conditions on the geometry of $X.$ 

\begin{theorem} 
Let $Z$ be a smooth projective 4-fold, $X,Y \subset Z$ smooth surfaces with G-positive normal bundles. Under one of the
following conditions $X$ and $Y$ meet.
\begin{enumerate}
\item $\kappa (X) = - \infty.$
\item $X$ is not minimal and every effective divisor $D$ in $Z$  has $\kappa (\sO_Z(D)) \geq 1.$
\end{enumerate} 
\end{theorem}

\begin{proof} (1) Choose a smooth rational curve $C \subset X$ with nef normal bundle $N_{C/X}.$ Since
$N_{X/Z}$ is ample, the normal bundle $N_{C/Z}$ is nef, hence the deformations of $C$ cover $Z$, in particular
some member of the family meets $Y.$ We conclude by (3.1). \\
(2) Choose a $(-1)-$curve $C \subset X.$ 
Using again the ampleness of $N_{X/Z} $ we conclude that either
$N_{C/Z}$ is nef or 
$$ N_{C/Z} = \sO(-1) \oplus \sO(a) \oplus \sO(b) $$
with $a,b > 0.$ 
In the first case we conclude as in (1). In the second we argue that the deformations of $C$ fill at least a divisor $D,$ see e.g. [Ko96,1.16]. 
In fact, assume the deformations cover only a surface $S.$ 
We consider a general member $C_t$ of the family of deformations of $C.$ We may assume that
$$ N_{C_t/Z} = \sO(-1) \oplus \sO(a') \oplus \sO(b') $$
with $a',b' > 0.$ Otherwise the normal bundle would be nef and the deformations of $C$ cover the whole $Z.$ 
Now choose a general smooth point $x \in S$ and a general $v \in T_{Z,x}$ which is normal to $S$. Then we find a section $s \in
H^0(N_{C_t/Z})$ such that $s(x) = v$ and therefore there is an infinitesimal deformation of $C_t$ along $v$. By non-obstructedness
this infinitesimal deformation extends to a deformation with positive-dimensional parameter space, so that we find deformations
of $C$ not contained in $S$, contradiction. \\
So the $C_t)$ fill a divisor $D$ (or the whole space, in which case we are done anyway). Since a multiple of $D$ moves by 
assumption, we conclude by
(3.2) that $D \cap Y \ne \emptyset.$ 
\end{proof}

We now treat Fano manifolds $Z.$ 
\begin{theorem} Let $Z$ be Fano 4-fold of index at least $2$, $X,Y \subset Z$ surfaces with $N_X$ ample and $N_Y$ G-positive. 
Then $X \cap Y \ne \emptyset.$ 
\end{theorem} 

\begin{proof} (1) We first treat the case $b_2(Z) = 1$ and give an argument which does not use classification. 
By Mella [Me99] (for index $2$, the index 3 case being settled by Fujita, see e.g. [IP99]), there is a smooth element 
$H \in \vert -K_Z \vert. $ Let $C $ be an irreducible component of $H \cap X.$ 
Then $C$ moves in an at least 1-dimensional family in the Fano 3-fold $H$. If the deformations of $C$ cover $H$, then some member of
the family meets $H \cap Y$, hence we conclude by (3.1). If the deformations of $C$ fill a divisor $D$ in $H$, then $D$ is ample in $H$, hence
$D \cap (H \cap Y) \ne \emptyset,$ and we conclude again by (3.1). \\
(2) In case $b_2(Z) \geq 2$ we need the classification of $Z$, see [Mu88,Mu89,IP99]. If $Z$ has index $3$, then $Z = \bP_2  \times \bP_2,$ hence homogenenous. 
If $Z$ has index $2$, either $Z$ is a product $\bP_1 \times W$ with $W = \bP_3$ or a del Pezzo 3-fold; hence we conclude by (5.1). 
Or $Z$ falls into of one 9 classes listed in [Mu88]. Then $Z$ is a divisor in a homogeneous manifold, a two-sheeted cover over a 
homogeneous manifold
or a $\bP_1-$bundle unless $Z$ is the blow-up $\phi$ of a 4-dimensional quadric $Q$ along a conic whose linear span is not contained in the quadric. 
In this case $Z$ has a quadric bundle structure over $\bP_2.$ Here we argue ad hoc as follows. We clearly have $\phi(X) \cap \phi(Y) \ne \emptyset.$
So if $X \cap Y = \emptyset,$ then both $X$ and $Y$ must meet $E$ (along a curve). Now $E = \bP_1 \times \bP_2,$ hence we can deform $X \cap E$
in $E$ to meet $Y \cap E.$ We conclude once more by (3.1). 
\end{proof} 

Addressing higher dimensions we first state

\begin{theorem} Let $Z$ be a del Pezzo manifold of dimension $n \geq 5$; $X$ and $Y$ submanifolds with $N_X$ ample and $N_Y$ G-positive such that
$\dim X + \dim Y \geq n.$ Then $X \cap Y \ne \emptyset.$ 
\end{theorem}

\begin{proof} Using Fujita's classification and the notation $-K_Z = (n-1)L,$ we are reduced to the following case: \\
$L^n = 1$ and $Z$ is a hypersurface of degree 6 in the weighted projective space $W = \bP(3,2,1, \ldots, 1).$ In this case we conclude by Proposition
5.6 below. \\
All other cases are 2-sheeted covers over projective spaces, hypersurfaces in homogeneous spaces or itself homogeneous. 
\end{proof} 

\begin{proposition} Let $Z \subset \bP(a_0, \ldots, a_n) $ be a smooth hypersurface in a weighted projective space $\bP(a_0, \ldots, a_n).$
Let  $X$ and $Y$ be submanifolds with $N_X$ ample and $N_Y$ G-positive such that
$\dim X + \dim Y \geq n.$ Then $X \cap Y \ne \emptyset.$ 
\end{proposition}

\begin{proof} We consider the projection $f: \bP_{n+1} \to  \bP(a_0, \ldots, a_n).$ By [Ba87,Prop.B]  any divisor $D_0$ in some irreducible component $X_0$
of $f^{-1}(X)$ moves inside a component $Z_0 $ of $f^{-1}(Z)$ containing $X_0$ such that the deformations $D_t$ cover $Z_0.$ Since $Z_0 \cap f^{-1}(Y) 
\ne \emptyset, $ there is some $t$ such that $D_t \cap f^{-1}(Y) \ne \emptyset. $ Thus the family $(f_*(D_t))$ deforms a divisor in $X$ to some 
$D_t$ which meets $Y.$ Henc $X \cap Y \ne \emptyset.$ 
\end{proof}   

We turn now to Fano manifolds $Z_n$ of index $n-2,$ so-called Mukai varieties. We will assume $n \geq 5$ and shall write $-K_Z = (n-2)H$; 
notice also the notion of the genus
of $Z$
$$ g = g(Z) = {{1} \over {2}}H^n + 1.$$
By [Mu88,89], $2 \leq g \leq 10.$

\begin{theorem}  Let $Z$ be a Fano of dimension $n \geq 5$ and index $n-2$.  
Let  $X$ and $Y$ be submanifolds with $N_X $ ample and $N_Y $ G-positive such that
$\dim X + \dim Y \geq n.$ Then $X \cap Y \ne \emptyset$  with the following possible exceptions.
\begin{enumerate}
\item $g = 5,$ \ $Z$ is the intersections of three quadrics in $\bP_{n+3},$ \ $ n = 2m$ and $\dim X = \dim Y = m.$
\item $ g= 7$, $5 \leq n \leq 8$ and $Z$ is a linear section of the $10-$dimensional rational-homogeneous manifold $SO_{10}(\bC)/P$ with $P$ maximal
parabolic.
\item $g = 8$, $5 \leq n \leq 6$ and $Z$ is a linear section  of the $8-$dimensional rational-homogeneous manifold $Sl_6(\bC)/P.$ 
\end{enumerate} 
\end{theorem}

\begin{proof} We shall use the classification due to Mukai ([Mu88,89], see also [IP99]). \\
If $b_2(Z) \geq 2,$ then $X = \bP_2 \times Q_3, \bP_3  \times \bP_3$ or
a hypersurface in $\bP_3 \times \bP_3,$ so we are done by (2.1) ($Q_n$ denotes the $n$-dimensional quadric). \\
So we shall assume $b_2(Z) = 1.$ 
In case $2 \leq g \leq 4,$ $Z$ is a degree 2 cover of $\bP_n$ 
resp. a hypersurface in the projective space or the quadric, hence our claim again holds by (2.1) and (4.12). If $g = 9,10$ again $Z$ is homogeneous or
a hypersurface in a homogeneous space, and we conclude. Thus it remains to treat the case $5 \leq g \leq 8.$ \\
In case $g = 5,$ we conclude from the Lefschetz hyperplane section theorem that $b_q(Z) = 1$ for all even $q \leq 2n$ with the exception $n = 2m$ and $q = m.$ 
Hence $X \cdot Y \ne \emptyset.$ \\
If $g = 6,$ then $Z$ is a degree 2 cover of $G(2,5),$ so we conclude by (4.11). \\
In the cases $g = 7,8,$ we can only treat the cases when $Z$ itself is homogeneous or a hyperplane of a homogeneous space. Thus only the listed cases
remain.

\end{proof}

\vskip 2cm \noindent
Thomas Peternell \\
Mathematisches Institut, Universit\" at Bayreuth \\
D-95440 Bayreuth, Germany \\
thomas.peternell@uni-bayreuth.de

\end{document}